\documentclass[a4paper,12pt]{amsart}

\usepackage{enumitem}

\usepackage{amsmath,amsthm,amsfonts,amssymb,mathtools,mathrsfs}

\usepackage{nicefrac}

\usepackage{tikz-cd}
\usetikzlibrary{babel}

\usepackage[utf8]{inputenc}

\usepackage[T1]{fontenc}
\usepackage{lmodern}

\usepackage{float}

\usepackage{subfigure}

\usepackage{xcolor}
\definecolor{colordelink}{rgb}{0,0,0.50}
\definecolor{colordecite}{rgb}{0,0.5,0}
\definecolor{colordeurl}{rgb}{0,0.41,0.5}

\usepackage[pagebackref=true]{hyperref}
\hypersetup{
    colorlinks=true,
    linkcolor=colordelink,
    citecolor=colordecite,
		urlcolor=colordeurl,
}

\usepackage[capitalize,nameinlink,noabbrev]{cleveref}

\usepackage[all]{hypcap}

\usepackage{todonotes}

\def\id{\operatorname{id}}

\def\ord{\operatorname{ord}}
\def\Re{\operatorname{Re}}
\def\Im{\operatorname{Im}}

\newcommand{\RR}{\mathbb{R}}

\newcommand{\CC}{\mathbb{C}}

\newcommand{\ZZ}{\mathbb{Z}}

\newcommand{\C}{\mathbb{C}}

\newcommand{\N}{\mathbb{N}}

\newcommand{\Ascr}{\mathscr{A}}

\newcommand{\abf}{\mathbf{a}}

\newcommand{\ebf}{\mathbf{e}}

\newcommand{\gbf}{\mathbf{g}}

\newcommand{\rbf}{\mathbf{r}}

\newcommand{\Sbb}{\mathbb{S}}

\newcommand{\mfrak}{\mathfrak{m}}

\newcommand{\Acal}{\mathcal{A}}

\newcommand{\Ccal}{\mathcal{C}}

\newcommand{\Ocal}{\mathcal{O}}

\makeatletter
\newcommand{\tpitchfork}{%
  \vbox{
    \baselineskip\z@skip
    \lineskip-.52ex
    \lineskiplimit\maxdimen
    \m@th
    \ialign{##\crcr\hidewidth\smash{$-$}\hidewidth\crcr$\pitchfork$\crcr}
  }%
}
\makeatother

\theoremstyle{plain}
\newtheorem{theorem}{Theorem}[section]
\newtheorem*{theorem*}{Theorem}
\newtheorem{lemma}[theorem]{Lemma}
\newtheorem*{lemma*}{Lemma}

\newtheorem*{corollary*}{Corollary}
\newtheorem{proposition}[theorem]{Proposition}
\newtheorem{question}[theorem]{Question}

\theoremstyle{definition}
\newtheorem{definition}[theorem]{Definition}
\newtheorem{proposition/definition}[theorem]{Proposition/Definition}

\newtheorem{example}[theorem]{Example}

\newtheorem{remark}[theorem]{Remark}

\Crefname{notation}{Notation}{Notations}

\theoremstyle{remark}

\begin{document}

\author{A. Fernández-Hernández and
R. Giménez Conejero}

\title[Complex plane curves up to diffeomorphism and rigidity]
{A note on complex plane curve singularities up to diffeomorphism and their rigidity}

\address{Departament de Matemàtiques,
Universitat de Val\`encia, Campus de Burjassot, 46100 Burjassot
SPAIN
\newline
Alfr\'ed R\'enyi Institute of Mathematics, Re\'altanoda utca 13-15,
H-1053 Budapest, 
HUNGARY
}
\email{alferher@alumni.uv.es
Roberto.Gimenez@uv.es
}



\begin{abstract} 
We prove that, if two germs of plane curves $(C,0)$ and $(C',0)$ with at least one singular branch are equivalent by a (real) smooth diffeomorphism, then $C$ is complex isomorphic to $C'$ or to $\overline{C'}$. A similar result was shown by Ephraim for irreducible hypersurfaces before, but his proof is not constructive. Indeed, we show that the complex isomorphism is given by the Taylor series of the diffeomorphism. We also prove an analogous result for the case of non-irreducible hypersurfaces containing an irreducible component of zero-dimensional isosingular locus. Moreover, we provide a general overview of the different classifications of plane curve singularities.
\end{abstract}

\maketitle

\thanks{The second author is partially supported by NKFIH Grant ``\'Elvonal (Frontier)'' KKP 144148.}
%


\section{Introduction}


A classical problem is the classification of plane curve singularities up to some equivalence. Among those, the most classical classification is the one given by ambient topological type or, simply, up to (ambient) homeomorphism. On the other side of the spectrum, both by rigidity and historically, one has the complex analytic classification. 

In this work we provide a self-contained proof of a classification that sits strictly in between these two classifications: classification up to (real) diffeomorphism of the ambient space. Given two germs of (not necessarily irreducible) plane curves $(C,0),(C',0)\subset(\CC^2,0)$ with equations $g(x,y)$ and $g'(x,y)$ and \textit{at least one singular branch}, we show in \cref{thm:main} the following.
\begin{theorem*}
$(C,0)$ and $(C',0)$ are diffeomorphic if, and only if, they are complex isomorphic or $C$ is  complex isomorphic to $\overline{C'}$, which is also a complex curve with equation $\overline{g'(\overline{x},\overline{y})}=0$.

 Furthermore, any diffeomorphism $\psi$ that takes one plane curve singularity to another, has a holomorphic or antiholomorphic Taylor series $T_\infty\psi$.
\end{theorem*}

This was previously done by Ephraim in \cite{Ephraim1973} for irreducible hypersurfaces, but with a non-constructive proof. Nevertheless, we present a stronger statement showing the rigidity of the diffeomorphism for curves. We also prove a partial result for non-irreducible hypersurfaces using Ephraim's ideas. It is important to know that the equivalent result for curves with only smooth branches is false, we show this in \cref{s:smoothbranches}.
\newline

In \cref{sec:class} we provide a concise review of the main different classifications we have, encompassing all the formal definitions and setting the notation. Furthermore, for some of them, some results in more dimensions are mentioned.

\cref{sec:power} is a short account on formal power series. Formal power series are also present in \cref{sec:formalparam}, where we study briefly what we call \textsl{formal parametrizations}.

The previous two sections are used in the proof of our main result, given in \cref{sec:main}. We also provide a detailed account on Ephraim's work in \cref{sec:ephraim} to show that these ideas are also useful to prove the non-irreducible setting in some cases.

We make some comments on the case of curves with only smooth branches in \cref{s:smoothbranches}, since both Ephraim's statement and the statement of our result is false in that case.

Finally, we offer extra insights in \cref{sec:conclude}.

\section{Different classifications}\label{sec:class}

While completing this work, we noticed that the singularity theory community lacks a small review about the main different classifications of plane curve singularities. Indeed, some of the results we state below are not widely known, despite their simplicity and importance. We include this short survey here, since we use the analytic classification, which is based on the topological classification, and give results about the smooth classification.

Unless otherwise stated, we restrict ourselves to the case of plane curve singularities. 

\begin{figure}[t]
	\centering
		\includegraphics[width=0.85\textwidth]{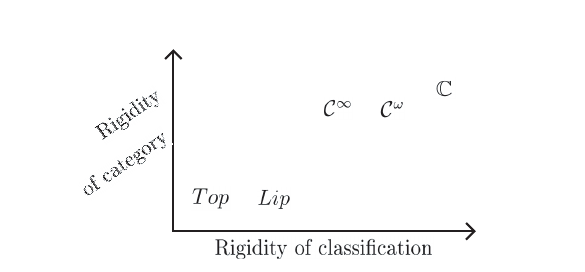}
	\caption{General overview of different classifications of plane curve singularities.}
	\label{fig:classifications}
\end{figure}

The general overview of the different classifications given by rigidity is provided in \cref{fig:classifications}.

\subsection{Topological}

The coarsest and most traditional classification is the topological classification. A comprehensive reference on this classification can be found, for example, in Wall's book \cite[Chapter 5]{Wall2004}, or also in \cite{Brieskorn1986,Pham1970}. See also a very detailed account in \cite{AlbertoTFG}. \\

\begin{definition}
We say that two curves $(C,0),(C',0)\subset(\CC^2,0)$ are \textit{(ambient) topologically equivalent} (also \textit{(ambient) homeomorphic}, or \textit{equisingular}) if there is a germ of homeomorphism $\psi:(\CC^2,0)\to(\CC^2,0)$ so that $\psi(C)=C'$. The class of equivalence is called \textit{topological class}.
\end{definition}

\begin{figure}[b]
	\centering
		\includegraphics[width=0.45\textwidth]{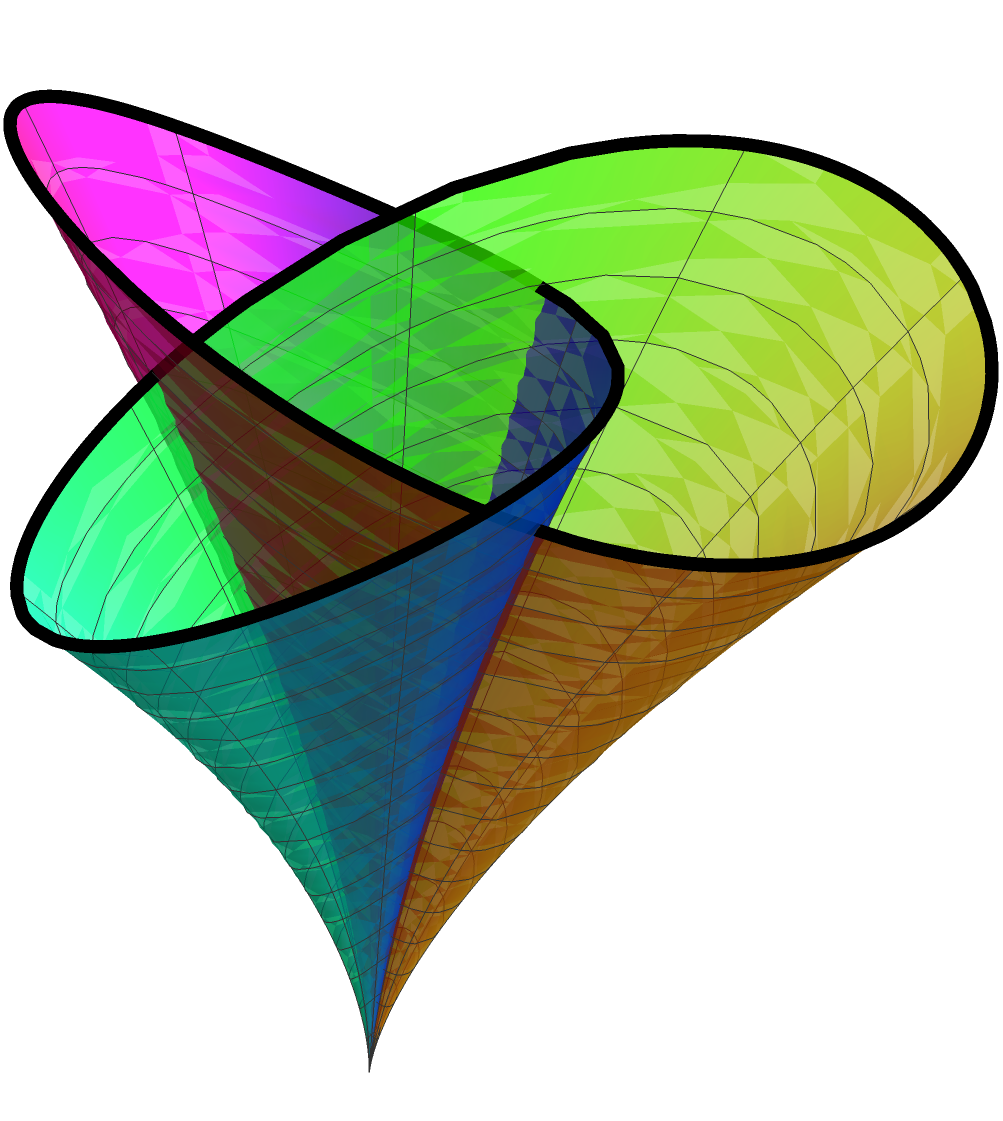}
      \includegraphics[width=0.50\textwidth]{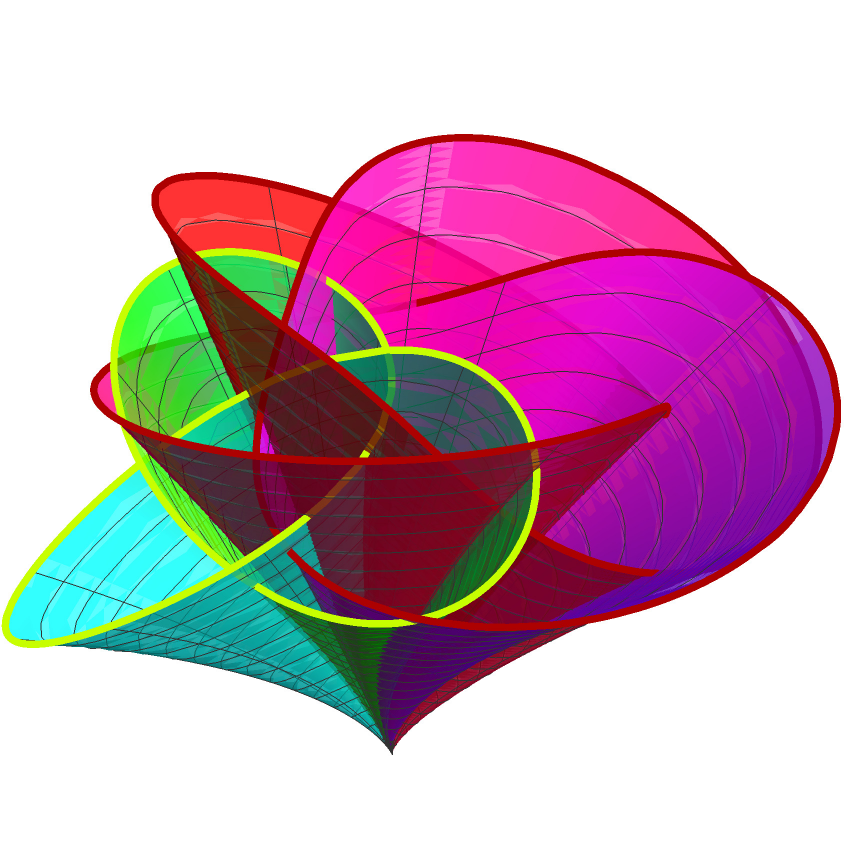}
	\caption{Representation of two plane curve singularities, one irreducible (left) and the other with two branches (right). Observe the conical structure over their link.}
	\label{fig:trefoilcurve}
\end{figure}

 In this case, it turns out that the topological class of a curve is completely determined by the isotopy class of its \textit{link}, namely the intersection of the curve with a sufficiently small sphere centered at the origin, because a curve is ambient homeomorphic to a cone over its link (see \cref{fig:trefoilcurve}). In this case, the \textsl{link} in the singularity theory sense is also a \textsl{link} in the topological sense, i.e., a disjoint union of knots in $\Sbb^3$ (one for each irreducible component of the curve). 

The origin of these ideas dates back to 1905, when Writinger started thinking about the possible groups that could arise from algebraic knots. He suggested this problem to his Ph.D. student Brauner, who gave the first steps towards the classification. They 
settled the key ideas to show that two branches with the same \textsl{Puiseux characteristic} have isotopic knots. The \textsl{Puiseux characteristic} and the \textsl{semigroup} of a branch are equivalent and contain the same information, so we define and mention the \textsl{semigroup} from now on.

\begin{definition}
Let $(C,0)$ be an irreducible plane curve (a branch) and assume $\gamma$ is an injective parametrization. The \textit{semigroup} of $C$ is the set of orders 
$$ \Gamma_C\coloneqq \big\{\ord \eta\circ\gamma: \eta\in\CC\left\{x,y\right\}\big\}.$$
\end{definition}

Brauner described the links via stereographic projections in his paper \cite{Brauner1928a}, in 1928. A clearer perspective, however, was given by Kähler in \cite{Kaehler1929}. Later, in 1932, Burau and Zariski noticed independently that the \textsl{Alexander polynomial} of an algebraic knot could be harnessed to extract the semigroup of its corresponding branch. This, therefore, showed that the semigroup is a complete topological invariant for branches.
More generally, a complete answer in the case of curves with several branches is the following theorem. A detailed historical account, with many references, is \cite[Theorem 1.1 of p. 730]{Kaehler2003}.
\begin{theorem}
    Two germs of analytic plane curves are topologically equivalent if, and only if, there exists a bijection between their branches that preserves both the semigroup of the branches and their pairwise intersection numbers.
\end{theorem}

\begin{remark}
Hypersurfaces in more dimensions are still classified by their link, since it is still true that the ambient topological type of a hypersurface singularity is that of the cone over its link (a general result in this direction is the \textsl{conic structure lemma} of Burghelea and Verona \cite[Lemma 3.2]{Burghelea1972}). The works of Saeki, Perron and King  also showed that topological equivalence of hypersurfaces coincides with \textsl{topological left-right equivalence} of their equations; \cite[Theorem 1]{Saeki1989}, \cite[Teórème 3']{Perron1985} and \cite[Theorem 3]{King1978}.
 A complete classification in terms of a mathematical object is far more complicated, and these kind of problems are still source of many other problems in different areas of singularity theory. Some examples of this are \cite{Le1973,Le1976}, the more recent \cite{Caubel2006,Nemethi2011,Nemethi2012} or the celebrated preprint \cite{Bobadilla2022}.
\end{remark}

\subsection{Bilipschitz}

The \textsl{Lipschitz} condition is slightly more restrictive than just continuity, but less than asking for analyticity or smoothness. Recall that a map $f$ between the metric spaces $(X,d_X)$ and $ (Y,d_Y)$ is \textit{Lipschitz} if there is a constant $c$ such that 
$$d_Y\big(f(p),f(q)\big)\leq c d_X(p,q).$$
Obviously, \textit{bilipschitz} maps are Lipschitz maps with a Lipschitz inverse and, since Lipschitz implies continuous, they are also homeomorphisms.

There are two natural metrics we can give to a complex set $(X,0)\subset (\CC^n,0)$ induced from the standard metric on $\CC^n$: the \textsl{outer metric} and the \textsl{inner metric}.

\begin{definition}
The \textit{outer metric} $d_{out}$ in $(X,0)$ is the one induced by the distance in $\CC^n$. The \textit{inner metric} $d_{inn}$ is given by the minimal length of (real) curves within $X$. They only depend on the complex analytic type of $X$ and not on the embedding.
\end{definition}

The first contributions on the classification of plane curves up to bilipschitz equivalence was given by Pham and Teissier in 1969 (see a translation in English in \cite{Pham2020}). Other contributions, including more dimensions, were also given by Mostowski in \cite{Mostowski1985}, Parusinski in \cite{Parusinski1994} and  Fernandes in \cite{Fernandes2003}. These results had some analytic, semi-algebraic or similar restrictions, but it was in 2013 that Neumann and Pichon in \cite{Neumann2013} that closed the classification for curves without extra restrictions. The complete result is the following.

\begin{theorem}
For two plane curves $(C,0),(C',0)$, the following are equivalent:
\begin{enumerate}[label=(\roman*),font=\itshape]
	\item there is a homeomorphism $\psi|:(C,0)\to(C',0)$ that is bilipschitz for the outer metric;
   \item there is a homeomorphism $\psi|:(C,0)\to(C',0)$, holomorphic except at zero, that is bilipschitz for the outer metric;
   \item there is a bilipschitz homeomorphism $\psi:(\CC^2,0)\to(\CC^2,0)$ such that $\psi(C)=C'$;
   \item $C$ is topologically equivalent to $C'$.
\end{enumerate}
\end{theorem}

For the inner metric it is very simple even for spatial curves.

\begin{proposition}
Any curve $(C,0)\subset(\CC^n,0)$ is metrically conical, i.e., bilipschitz to the metric cone on its link, for the inner metric.
\end{proposition}

\begin{remark}
Lipschitz geometry has experienced an explosion of activity in the recent years in many directions. For instance, the global case is also studied, see \cite{Fernandes2020,Targino2023,Sampaio2023}, but many other aspects of Lipschitz geometry are also being studied, for example \cite{Bobadilla2017,Birbrair2018}.
\end{remark}

\subsection{Real analytic and smooth}

Regarding rigidity of \textsl{isomorphisms}, the next ones among the typical ones are equivalences by smooth (real) diffeomorphism and real analytic diffeomorphism. This is the object of study of this paper.

\begin{definition}
We say that two plane curve singularities $(C,0),(C',0)$ are \textit{(ambient) diffeomorphic}, \textit{smooth equivalent} or $\Ccal^\infty$\textit{-isomorphic}
 if there is a (real) diffeomorphism $\psi:(\CC^2,0)\to(\CC^2,0)$ such that $\psi(C)=C'$.

Analogously, we define the notion of \textit{real analytic diffeomorphic} or $\Ccal^\omega$\textit{-isomorphic} curves when the diffeomorphism is real analytic.
\end{definition}

These two of equivalences were first addressed by Ephraim in his Ph.D. thesis in 1971, under the supervision of Laufer and Gunning, later published as an article in \cite{Ephraim1973}. He showed that the smooth and real analytic equivalences are \textsl{equivalent} themselves for real analytic sets, in full generality. In other words:

\begin{proposition}
Let $V$ and $W$ be germs of real analytic sets. They are $\Ccal^\infty$-isomorphic if, and only if, they are $\Ccal^\omega$-isomorphic.
\end{proposition}

Now, observe that if $V$ is a complex set, then $\overline{V}$ is also a complex set, given by the same equation but with conjugated coefficients. Ephraim also showed the following.

\begin{theorem}\label{thm:introephra}
Let $(V,0)$ be an irreducible germ of hypersurface. Then, if $V$ is $\Ccal^\infty$-isomorphic to a germ of a complex set $W$, $V$ is $\CC$-isomorphic to $W$ or $\overline{W}$.
\end{theorem}

We review his proof in \cref{sec:ephraim}, and use the same ideas to show the following.

\begin{proposition}
\cref{thm:introephra} holds also for non-irreducible hypersurfaces if $V$ contains a component that cannot be written as a product of varieties of lower dimension.
\end{proposition}

Finally, we give a stronger statement for the case of plane curve singularities in \cref{sec:main} with at least one singular branch. The case of smooth branches is covered in \cref{s:smoothbranches}, showing that the equivalent statement of the previous \cref{thm:introephra} and the following theorem fail in different ways.

\begin{theorem}
If $(C,0)$ is $\Ccal^\infty$-isomorphic to $(C',0)$ by a diffeomorphism $\psi$, its Taylor series $T_\infty\psi$ is holomorphic or antiholomorphic.
\end{theorem}

\begin{remark}
It is surprising how little these direction is known and how little it has been researched. There are few instances of works focusing on the smooth aspect of complex varieties. For example, Ephraim showed that one can write any irreducible complex variety of any codimension as a product of \textit{indecomposable} varieties and used this structure theorem to show when two complex sets are $\Ccal^\infty$-isomorphic in \cite{Ephraim1976}. He also generalized his structure theorem for any reduced variety in \cite{Ephraim1978}.
\end{remark}

\subsection{Complex analytic}

The most restrictive classification is, naturally, the complex analytic classification.

\begin{definition}
Two plane curves $(C,0),(C',0)$ are \textit{biholomorphic}, $\CC$\textit{-isomorphic} or \textit{analytic equivalent} if there is a biholomorphism $\phi:(\CC^2,0)\to(\CC^2,0)$ such that $\phi(C)=C'$.
\end{definition}

 While the mathematical community has held an interest in the classification of plane curves for a considerable period, a complete classification was not attained until 2020. The first complete classification was given for branches in 2011 by Hefez and Hernandes (Marcelo Escudeiro) in \cite{Hefez2011,Hefez2021}. For two branches by Hefez, Hernandes (Marcelo Escudeiro) and Hernandes (Maria Elenice Rodrigues) in \cite{Hefez2015} in 2015. And, finally, for non-irreducible curves in a preprint version of 2021 by Hernandes (Marcelo Escudeiro) and Hernandes (Maria Elenice Rodrigues) in \cite{Hernandes2020}. 

A huge number of authors have contributed to this problem before these works were made. To mention some of them together with some of their work: Ebey gave some normal forms in 1965 in \cite{Ebey1965}; Zariski made important contributions, in particular to the moduli space of some class of curves in 1973 in \cite{Zariski1973}; Delorne, Granger, Laudal and Pfister, and Greuel and Pfister also contributed to the study of the moduli space in, respectively, \cite{Delorme1978,Granger1979,Laudal1988,Greuel1994}; Bruce and Gaffney classified simple irreducible curves in 1982 in \cite{Bruce1982}; Kolgushkin and Sadykov, and Genzmer and Paul gave normal forms in some cases in, respectively, \cite{Kolgushkin2001,Genzmer2011}; and many more.

\begin{definition}
Let $(C,0)$ be a plane branch parametrized by $\gamma(t)=\big(\gamma_1(t),\gamma_2(t)\big)$. The \textit{set of orders of Kähler differentials} of $C$ is 
$$ \Lambda\coloneqq \big\{\ord_t \big(\eta\circ\gamma(t) \gamma_1'(t)+\xi\circ\gamma(t)\gamma'_2(t)\big)+1:\eta,\xi\in\Ocal_2\big\}\setminus\left\{1\right\}.$$
It is a $\Gamma_C$-module, i.e., $\Gamma_C+\Lambda\subseteq\Lambda$,
 and it is $\Ascr$-invariant. Furthermore, if $\Lambda\setminus\Gamma\neq\varnothing$, the \textit{Zariski invariant} of $C$ is
$$ \lambda\coloneqq \min(\Lambda\setminus\Gamma)-v_0,$$
where $v_0$ is the multiplicity of the curve.
\end{definition}

When it comes to study irreducible plane curves under $\CC$-isomorphism, it is a well-known fact that the isomorphism class of the branches coincides exactly with the $\Ascr$-equivalence class of their parametrizations. 

\begin{theorem}[Normal Forms Theorem, see {\cite[Theorem 2.1]{Hefez2011}}]\label{thm: normalform}
Let $\gamma$ be an injective parametrization of a singular plane branch with semigroup $\Gamma=\left\langle v_0,\dots,v_\ell\right\rangle$. Then, either $\gamma$ is $\Ascr$-equivalent to the monomial parametrization $(t^{v_0},t^{v_1})$ or to the parametrization
\begin{equation}
\left(t^{v_0},t^{v_1}+t^\lambda+\sum_{i>\lambda, i\notin\Lambda-v_0}a_it^i\right),
\label{eq:Normalparametrization}
\end{equation}
where $\lambda$ is its Zariski invariant and $\Lambda$ is the set of orders of Kähler differentials of the branch.

Moreover, if $\gamma$ and $\gamma'$, with coefficients $a_i'$ instead of $a_i$, are parametrizations as in \cref{eq:Normalparametrization} representing two plane branches with same semigroup and same set of values of differentials, then $\gamma \sim_\Ascr \gamma'$ if and only if there is $r\in\CC^*$ such that $r^{\lambda -v_1}=1$ and $a_i=r^{i-v_1}a_i'$, for all $i$.
\end{theorem}

 This classification is far more subtle, since the moduli space of a curve can be non-discrete, i.e., a family $\{C_\varepsilon\}_\varepsilon$ varying analytically in $t$ can yield different analytic classes for each $t$ (an example of this phenomenon is given in \cref{ex:moduli}). In contrast, the topological class of a plane branch parametrized by 
\begin{equation*}
    \gamma (t)=(t^{v_0}, \sum_{r> m} a_r t^r)
\end{equation*}
is solely determined by its exponents $r$ with $a_r\neq0$ (specifically those in the \textsl{Puiseux characteristic}) instead of the coefficients $a_r$ themselves.


\subsection{Other classifications}

There are some categories that we have not mentioned, e.g., the one regarding $\Ccal^1$ homeomorphisms. As far as we know, not much about this has been studied. For example, it is known that multiplicity of hypersurfaces is a $\Ccal^1$-invariant (see \cite{Ephraim1976a}). Examples seem to suggest that the $\Ccal^1$-isomorphism classification may coincide with the topological one (see \cref{ex:c1notan,qu:c1} below).

There are other equivalences, specially in the real analytic or the global setting. See for example the \textsl{blow-spherical equivalence} in \cite{Birbrair2017,Sampaio2022}, and the blow-analytic equivalence in \cite{Kuo1985,Kobayashi1998,Fichou2011,Valle2016}.

\section{Preliminaries on formal power series}\label{sec:power}

\subsection{Formal power series}\label{sec:Faa}
We will use fluently some fundamental facts about formal power series so, for the sake of clarity, we state them here.
\newline

\begin{proposition}\label{prop:faadibruno}
Given the formal power series $\mathbf{a}(t)=\sum_{n=0}^\infty a_n t^n$ and $\mathbf{b}(t)=\sum_{n=0}^\infty b_n t^n$, where $b_0=0$, one has that the composition is
$$ \mathbf{c}(t)\coloneqq\mathbf{a}\big(\mathbf{b}(t)\big)=\sum_{n=1}^\infty c_n t^n,$$
where $c_0=a_0$ and
\begin{equation}\label{eq:coefficients}
 c_n = \sum_{k=1}^n a_k\sum_{\mathbf{i}\in P_k(n)}b_{i_1}\cdots b_{i_k},
\end{equation}
for the partitions
$$P_k(n)\coloneqq\left\{(i_1,\dots,i_k): i_1+\cdots+i_k=n\right\}.$$
\end{proposition}

This is commonly known as \textit{Faà di Bruno's formula}. With it, knowing two of the series $\mathbf{a},\mathbf{b}$ or $\mathbf{c}$ is sometimes enough to compute the other by means of \cref{eq:coefficients}.

\begin{definition}
We say that the \textit{order} of $\mathbf{a}$ is $v_0$ if $a_0=\cdots=a_{v_0-1}=0$ and $a_{v_0}\neq0$. 

Furthermore, a \textit{formal change of variable} on the series $\mathbf{a}$ is precomposing with a series $\mathbf{b}$ of order one.
\end{definition}


\begin{lemma}\label{lem:changetv0}
    If $\abf\in \C [[t]]$ is a formal power series of order $n$, then there exists a formal power series of order one $\rbf$ such that $\rbf^n=\abf$. Furthermore, if $\abf$ has order $v_0$, it admits a formal change of variable so that the resulting series is precisely $t^{v_0}$.
\end{lemma}
\begin{proof} Since $\abf$ has order $n\geq 1$, one can write $\abf=ax^n(1+\ebf)$ for $a\neq 0$ and $\ebf$ a formal power series with constant term $0$. Let $r\in \C$ be such that $b^n=a$, and expand $\ebf=\sum_{s=1}^{+\infty} e_st^s$. Notice that we have an expansion for $(1+\ebf)^{1/n}$ given by 
    \begin{equation}\label{eq:alberto}
        (1+\ebf)^{n^{-1}}=1+n^{-1}\ebf+\dfrac{n^{-1}(n^{-1}-1)}{2}\ebf^2+ \dots  
    \end{equation}

    For a fixed $s\in \N$, only a finite number of terms of \cref{eq:alberto} contribute to the coefficient of $t^s$. Thus, the previous expansion defines a power series $\gbf(t)$. Hence the series $\rbf=bt\gbf(t)$ has clearly order $1$ since the initial term of $\gbf$ is $1$ and $b\neq 0$. It satisfies that $\rbf^n=at^n(1+\ebf)=\abf$. 
   
   Finally, as an immediate consequence, any power series of order $v_0$ admits a formal change of variable so that the resulting series is precisely $t^{v_0}$ by the \textsl{Lagrange inversion formula} (see for example \cite[Proposition 5.4.1 and Theorem 5.4.2]{Stanley1999}). 
\end{proof}

\begin{remark}\label{thm:fields}
The existence of roots of power series may be easily derived from this result. More precisely, if $\abf\in \CC[[t]]$ has order $mn$, then $\abf$ may be written as $\abf=\rbf^{mn}$ for some formal power series $\rbf$ of order one. Then $\rbf^n$ is an $m$-th root of $\abf$.
\end{remark}

\section{Formal parametrizations}\label{sec:formalparam}

We need to prove some results on what we call \textsl{formal parametrizations}. 

\begin{definition}\label{def:formalparam}
A \textit{formal parametrization} $\gamma_\infty$ of a curve $(C,0)\subset(\CC^2,0)$ given by a complex equation $g$ is a pair of formal power series
$$ \gamma_\infty(t)=\left(\sum_{i=1}^\infty\alpha_i t^i,\sum_{i=1}^\infty\alpha_i' t^i\right) $$
such that the formal composition $g\circ\gamma_\infty$ is zero.
\end{definition}

\begin{example}\label{ex:y=0}
Formal parametrizations are not necessarily analytic. For example, consider any analytic parametrization $\gamma(t)$ of a curve $C$ and make a change $t$ for any badly-behaved formal series. The easiest, and most disappointing, example happens for $\gamma(t)=(t,0)$, the curve $\left\{y=0\right\}$ and the series $\sum n!t^n$.
\end{example}

For a given formal parametrization, if we try to \textsl{undo} the change shown in \cref{ex:y=0} by a formal change of variables (which we can do by \cref{thm:fields}), we recover something analytic by the lemma below. This is particularly relevant when we do not know if the formal parametrization comes from a change such as the one in \cref{ex:y=0}.

\begin{lemma}[cf. {\cite[Corollary 5.1.8 (3)]{Jong2000}}]\label{lem:closedfield}
A formal parametrization $\gamma_\infty(t)$ of the form
$$\gamma_\infty(t)=\left(t^m,\sum_{i=1}^\infty\alpha_i' t^i\right)$$
is analytic. The same holds if, more generally, the first formal power series is analytic.
\end{lemma}


It is a well-known fact that every injective parametrization $\gamma:(\C,0)\rightarrow (\C^2,0)$ is $\Ascr$-finite in Mather's sense. Therefore, for large enough $k$, $\gamma$ and its Taylor polynomial of order $k$ parametrize $\CC$-isomorphic branches  by Mather finite determinacy theorem (see \cite[Chapter 6]{Mond2020}).

\begin{theorem}\label{thm:kjet}
Let $\gamma_\infty(t)$ be a formal parametrization of the curve $C$ of minimal order (the multiplicity of $C$). Assume it has the form
$$\gamma_\infty(t)=\left(\sum_{i=1}^\infty\alpha_i t^i,\sum_{i=1}^\infty\alpha_i' t^i\right).$$
Then, for a big enough $k$,
$$\gamma_k(t)\coloneqq\left(\sum_{i=1}^k\alpha_i t^i,\sum_{i=1}^k\alpha_i' t^i\right)$$
parametrizes a curve isomorphic to $C$.
\end{theorem}
\begin{proof}
Assume that the first formal power series has lower order than the second one, and that this order is $m$, i.e., $\alpha_1=\dots=\alpha_{m-1}=0$ and $ \alpha_m\neq0$.
Then, we make a change of variables so that
$$ \tilde{t}^m=\sum_{i=m}^\infty\alpha_i t^i.$$
In other words, we have to take an $m$-root of the power series (i.e., a solution of $\sum_{i=m}^\infty\alpha_i t^i-y^m=0$, which exist by \cref{thm:fields}), say $\rbf(t)$, and its inverse by composition $\rbf^{\left\langle -1\right\rangle}$, as it exists since $\rbf$ has order 1 (by the Lagrange inversion formula, see for example \cite[Proposition 5.4.1 and Theorem 5.4.2]{Stanley1999}). As we want $\tilde{t}=\rbf(t)$, we have to perform the change $t=\rbf^{\left\langle -1\right\rangle}\big(\tilde{t}\big)$. 

After this, as the new formal parametrization is also a formal parametrization of $C$ with first entry $\tilde{t}^m$, we have a complex parametrization $\gamma(\tilde{t})$ of $C$, by \cref{lem:closedfield}. Notice that both $\rbf$ and its inverse $\rbf^{\left\langle -1\right\rangle}$ can be computed recursively so, for two series that are equal up to a high order, the $m$-roots and their inverses are equal up to a high degree as well (up to choice of $m$-roots). This is true, in particular, if we take 
$$ \sum_{i=m}^k\alpha_i t^i \quad\text{ and }\quad \sum_{i=m}^\infty\alpha_i t^i.$$
Hence if, for $k\gg0$, we consider
$$\gamma_k(t)\coloneqq\left(\sum_{i=1}^k\alpha_i t^i,\sum_{i=1}^k\alpha_i' t^i\right) $$ and proceed similarly with the polynomial $\sum_{i=1}^k\alpha_i t^i$, we reach another complex parametrization that is equal, up to some degree, to $\gamma(\tilde{t})$. In particular, since $\gamma_k$ is $\Ascr$-finite by the assumption on the order of $\gamma$, by Mather finite determinacy theorem (see \cite[Chapter 6]{Mond2020}), the parametrizations are $\Ascr$-equivalent, so the curves $C$ and the curve that is parametrized by $\gamma_k$ are isomorphic.
\end{proof}

\section{Main result}\label{sec:main}
   
   Let us consider two irreducible plane curve singularities, $(C,0)$ and $(C',0)$ in $(\CC^2,0)$, with equations $g$ and $g'$,  that are $\Ccal^\infty$-isomorphic by the diffeomorphism $\psi(x_R,x_I,y_R,y_I)$. By composing with convenient $\CC$-isomorphisms, we can assume that the curves admit a normal form parametrization (see \cref{thm: normalform}):
\begin{equation*}
\begin{aligned}
\gamma(t)&=\left(t^{v_0},t^{v_1}+t^\lambda+\sum_{s>\lambda, \:s\notin\Lambda-v_0}\omega_st^s\right), \text{ and }\\
\gamma'(t)&=\left(t^{v_0},t^{v_1}+t^{\lambda'}+\sum_{s>\lambda', \:s\notin\Lambda'-v_0}\omega_s't^s\right).
\end{aligned}
\end{equation*} 
To simplify notation, we set $\omega_{v_1}=\omega_{\lambda}=\omega_{\lambda'}\coloneqq 1$ and we write $\gamma(t)=\left(t^{v_0},\sum \omega_st^s\right)$ and $\gamma'(t)=\left(t^{v_0},\sum \omega_s't^s\right)$ from now on.

If we consider a real line in $(\CC,0)$ parametrized by $te^{i\theta}$ for $t\in(\RR,0)$ and fixed $\theta\in\RR$, we have that
$$ g'\circ \psi\circ\gamma\big(te^{i\theta})=0. $$
Hence, the same holds true for the composition of Taylor series
$$ g'\circ T_\infty\psi\circ\gamma(te^{i\theta})=0. $$
Therefore, as $T_\infty\psi\circ\gamma(te^{i\theta})$ is, for every $\theta\in\RR$, a formal power series on $t$ that satisfies the equation $g'$, $T_\infty\psi\circ\gamma(te^{i\theta})$ is a formal parametrization of $C'$ as defined in \cref{def:formalparam}. 

We will study this formal parametrization to see that, indeed, $T_\infty\psi$ has to be very simple provided, as it is with our assumptions, that $C$ and $C'$ are in normal forms.

Let us establish the notation. The Taylor series of $\psi$ is the formal power series
\begin{equation*}
\begin{aligned}
T_\infty\psi (x_R,x_I,y_R, y_I)=&\Bigg(\sum_{i,j,k,\ell}\alpha^R_{ijk\ell}x_R^ix_I^jy_R^ky_I^\ell,\sum_{i,j,k,\ell}\alpha^I_{ijk\ell}x_R^ix_I^jy_R^ky_I^\ell,\\
                                                          &\quad\sum_{i,j,k,\ell}\beta^R_{ijk\ell}x_R^ix_I^jy_R^ky_I^\ell,\sum_{i,j,k,\ell}\beta^I_{ijk\ell}x_R^ix_I^jy_R^ky_I^\ell\Bigg),
\end{aligned}
\end{equation*} 
which, writing $\alpha_\bullet\coloneqq\alpha_\bullet^R+i\alpha_\bullet^I$ and $\beta_\bullet\coloneqq\beta_\bullet^R+i\beta_\bullet^I$, can be written as
\begin{equation*}
T_\infty\psi (x_R,x_I,y_R, y_I)=\Bigg(\sum_{i,j,k,\ell}\alpha_{ijk\ell}x_R^ix_I^jy_R^ky_I^\ell,\sum_{i,j,k,\ell}\beta_{ijk\ell}x_R^ix_I^jy_R^ky_I^\ell\Bigg).
\end{equation*} 
   Without any further information, the formal parametrization has the form
   \begin{equation*}
\begin{aligned}
&T_\infty\psi\circ\gamma(te^{i\theta})= \\
&\resizebox{0.97\hsize}{!}{%
$\Bigg( {\displaystyle\sum_{i,j,k,\ell}} \alpha_{ijk\ell} \cos^i(\theta v_0) \sin^j(\theta v_0)t^{v_0(i+j)}\left(\textstyle \sum_s \Re(\omega_s e^{i\theta s})t^s \right)^k\left(\textstyle \sum_s \Im(\omega_s e^{i\theta s})t^s \right)^\ell$}%
;\\
 & \resizebox{0.95\hsize}{!}{%
${\displaystyle\sum_{i,j,k,\ell}} \beta_{ijk\ell} \cos^i(\theta v_0) \sin^j(\theta v_0) t^{v_0(i+j)}\left(\textstyle \sum_s \Re(\omega_s e^{i\theta s})t^s \right)^k\left(\textstyle \sum_s \Im(\omega_s e^{i\theta s})t^s \right)^\ell\Bigg)$}%
.
\end{aligned}
\end{equation*}
 \subsection{Simplification to linear terms} \label{subsec:linear}
   By \cref{lem:changetv0}, we can find a formal change of variable such that the first entry of the formal parametrization $T_\infty\psi\circ\gamma(te^{i\theta})$ is $t^{v_0}$. If we proceed with this change of variable, by \cref{lem:closedfield}, we necessarily recover a normal form parametrization $\gamma'(t)=(t^{v_0},\sum \omega'_s t^s)$ (or any other normal form equivalent to this one). 
   
   As we pointed out in \cref{sec:Faa}, if we know the result of the change of variables we obtain some information of the change of variables itself by means of the equations of Faà di Bruno's formula, \cref{prop:faadibruno}. Following the notation of \cref{sec:Faa}, we have
   \begin{equation*}
\begin{aligned}
\mathbf{a}(t)&=\resizebox{0.9\hsize}{!}{%
$\sum \alpha_{ijk\ell} \cos^i(\theta v_0) \sin^j(\theta v_0)t^{v_0(i+j)}\left(\textstyle \sum_s \Re(\omega_s e^{i\theta s})t^s \right)^k\left(\textstyle \sum_s \Im(\omega_s e^{i\theta s})t^s \right)^\ell$}%
\\
\mathbf{c}(t)&=t^{v_0}
\end{aligned}
\end{equation*}
   and we want to find $\mathbf{b}(t)$, the change of variables. As it is the same change we perform on the second entry of the formal parametrization and it has to yield $\sum \omega'_s t^s$, we will obtain information about the first series $\mathbf{a}(t)$.

   The first equations of Faà di Bruno's formula are:
\begin{equation}\label{eq:FDB1}
\begin{aligned}
c_1=0&=0,\\
\vdots&\quad\vdots\\
c_{v_0-1}=0&=0,\\
 c_{v_0}=1&= a_{v_0} b_1^{v_0}.
\end{aligned}
\end{equation}
Assume that $\mathbf{a}(t)$ has a term of order bigger than $v_0$. This happens if, and only if, there is some $\alpha_{ijk\ell}\neq0$ with $(i,j,k,\ell)\neq(1,0,0,0),(0,1,0,0)$. Then, if $M$ is the minimal from such orders, the equations continue with
\begin{equation*}
\begin{aligned}
c_{v_0+1}=0&=\kappa a_{v_0}  b_1^{v_0-1}b_2+ a_{v_0+1}b_1^{v_0+1},
\end{aligned}
\end{equation*}
for some irrelevant constant $\kappa$, giving $b_2=0$ if $v_0+1<M$ (so $a_{v_0+1}=0$). For the same reason, $b_3=\dots=b_{M-v_0}=0$. Finally, the equation
\begin{equation*}
\begin{aligned}
 c_{M}=0&= \kappa' a_{v_0}  b_1^{v_0-1}b_{M-v_0+1}+a_M b_1^{M}
\end{aligned}
\end{equation*}
   gives $b_{M-v_0+1}\neq0$. This term in the change of variable will reflect into restrictions on $\mathbf{a}(t)$ when we proceed with the change of variable of the second entry of the formal parametrization $T_\infty\psi\circ\gamma(te^{i\theta})$.  Indeed, by the expression of $\mathbf{a}(t)$, 
   $$M=v_0(i_0+j_0)+v_1(k_0+\ell_0)$$
   for some $i_0,j_0,k_0,\ell_0$, so $b_{v_0(i_0+j_0-1)+v_1(k_0+\ell_0)+1}\neq0$.
   
   When we consider the change of variable given by $\mathbf{b}(t)$ on the second entry, the resulting series contains a term, say $T$, of order 
   $$v_1-1+(M-v_0+1)=v_0(i_0+j_0-1)+v_1(k_0+\ell_0+1)$$ coming from a product of $v_1-1$ times $b_1$ and one time $b_{M-v_0+1}$ when computing $\mathbf{b}^{v_1}$. Indeed, $T$ has the minimum order one can get after the change of variables if one ignores all the terms coming uniquely from the linear part $b_1$ of $\mathbf{b}(t)$. We will use this fact shortly.
   
   Recall that the resulting change of variable yields
   $$\gamma'(t)=\left(t^{v_0},t^{v_1}+t^{\lambda'}+\sum_{s>\lambda', \:s\notin\Lambda'-v_0}\omega_s't^s\right).$$
   Hence, as $v_0(i_0+j_0-1)+v_1(k_0+\ell_0+1)\in\Gamma'-v_0\subset\Lambda'-v_0$, the resulting term of this order has to be zero (for every $\theta$) since it is not $v_1,\lambda$ nor any $s\notin\Lambda'-v_0$. Furthermore, there are no cancellations because $T$ has minimum order, so $T=0$ for every $\theta$. This is a contradiction. Therefore, $\alpha_{ijk\ell}=0$ if $(i,j,k,\ell)\neq(1,0,0,0),(0,1,0,0)$.
   
   Hence, the formal parametrization has the form
   \begin{equation*}
\begin{aligned}
&T_\infty\psi\circ\gamma(te^{i\theta})= \\
&\bigg( \big(\alpha_{1000} \cos(\theta v_0) +\alpha_{0100} \sin(\theta v_0) \big)t^{v_0}%
;\\
 & \resizebox{0.95\hsize}{!}{%
${\displaystyle\sum_{i,j,k,\ell}} \beta_{ijk\ell} \cos^i(\theta v_0) \sin^j(\theta v_0) t^{v_0(i+j)}\left(\textstyle \sum_s \Re(\omega_s e^{i\theta s})t^s \right)^k\left(\textstyle \sum_s \Im(\omega_s e^{i\theta s})t^s \right)^\ell\bigg)$}%
.
\end{aligned}
\end{equation*}   
It is obvious that the change of variables $\mathbf{b}(t)$ is just a constant times $t$. After the change of variables we recover $\gamma'(t)$, so $\beta_{ijk\ell}=0$ if $i+j\geq1$ or $k+\ell\geq2$, otherwise we would have a term of order $v_0(i+j)+v_1(k+\ell)\in\Lambda'-v_0$. Therefore, the parametrization has the form
   \begin{equation*}
\begin{aligned}
T_\infty\psi\circ\gamma(te^{i\theta})= &\bigg( \big(\alpha_{1000} \cos(\theta v_0) +\alpha_{0100} \sin(\theta v_0) \big)t^{v_0};\\
 &\beta_{0010} \left(\textstyle \sum_s \Re(\omega_s e^{i\theta s})t^s \right)+\beta_{0001}\left(\textstyle \sum_s \Im(\omega_s e^{i\theta s})t^s \right)\bigg)%
\end{aligned}
\end{equation*} 
and
$$ T_\infty\psi (x_R,x_I,y_R, y_I)=\left(\alpha_{1000}x_R+\alpha_{0100}x_I,\beta_{0010}y_R+\beta_{0001}y_I\right). $$

\subsection{Holomorphy} We are going to prove that $T_\infty\psi$ is either holomorphic or antiholomorphic (i.e., holomorphic after conjugation on both variables) by studying further the change of variable we have shown above. To simplify notation, we set
   \begin{equation*}
\begin{aligned}
\alpha\coloneqq&\alpha_{1000}\\
\alpha'\coloneqq&\alpha_{0100}\\
\beta\coloneqq&\beta_{0010}\\
\beta'\coloneqq&\beta_{0001}.
\end{aligned}
\end{equation*} 

\begin{lemma}\label{lem:orthogonal}
The coefficients $\alpha,\alpha',\beta,\beta'$ are not zero. Furthermore, $\alpha$ is orthogonal to $\alpha'$ and $\beta$ is orthogonal to $\beta'$.
\end{lemma}
\begin{proof}
First of all, as $\psi$ is a diffeomorphism,  its differential
\[d\psi_0=\begin{bmatrix} 
 \alpha^R   &  \alpha'^R  &  0 &0  \\
 \alpha^I   &  \alpha'^I  &  0 & 0 \\
 0   &  0  &  \beta^R   &  \beta'^R  \\
 0   &  0  &  \beta^I   &  \beta'^I  
\end{bmatrix},
\]
has maximum rank, so 
$\alpha, \alpha', \beta, \beta'$ are non-zero.

From \cref{eq:FDB1} we see that
$$ b_1=\frac{1}{a_{v_0}^{\nicefrac{1}{v_0}}}=\frac{1}{\big(\alpha \cos(\theta v_0) +\alpha' \sin(\theta v_0) \big)^{\nicefrac{1}{v_0}}}.$$
We do not care at the moment about the lack of precision about the $v_0$-th root.

When we apply the change of variable given by $\mathbf{b}(t)=b_1 t$ on the second entry of $T_\infty\psi\circ\gamma(te^{i\theta})$, we see that, necessarily, the coefficient of $t^{v_1}$ is
$$ 1= \frac{\beta\cos(\theta v_1)+\beta'\sin(\theta v_1)}{\big(\alpha \cos(\theta v_0) +\alpha' \sin(\theta v_0) \big)^{\nicefrac{v_1}{v_0}}},$$
so
\begin{equation}
\big(\alpha \cos(\theta v_0) +\alpha' \sin(\theta v_0) \big)^{v_1}=\big(\beta\cos(\theta v_1)+\beta'\sin(\theta v_1)\big)^{v_0}.
\label{eq:linear1}
\end{equation}
Setting $\theta=\pm\frac{\pi}{2v_0}$ in \cref{eq:linear1} we obtain
$$ (\pm1)^{v_1}\alpha'^{v_1}=\left(\beta\cos\left(\frac{\pi v_1}{2v_0}\right)\pm\beta'\sin\left(\frac{\pi v_1}{2v_0}\right)\right)^{v_0}.$$
Taking modules, we see that $\beta$ and $\beta'$ have to be orthogonal. Likewise, with $\theta=\pm\frac{\pi}{2v_1}$, one deduces that $\alpha$ and $\alpha'$ are orthogonal as well.
\end{proof}

We need a small technical lemma.
\begin{lemma}\label{lem:orthogonalcomplex}
For a pair of orthogonal complex numbers $z,z'\in\CC\setminus0$ and $\rho\notin\frac{\pi}{2}\ZZ$, if
$$ \frac{\left|z\cos\rho+z'\sin\rho\right|}{\left|z\sin\rho+z'\cos\rho\right|}=\frac{\left|z\right|}{\left|z'\right|} $$
then $\left|z\right|=\left|z'\right|$.
\end{lemma}
\begin{proof}
After a change of coordinates, we can assume that $z=a$ and $z'=ib'$ for $a,b'\in\RR$. Then, taking squares in both sides,
$$ \frac{a^2\cos^2\rho+b'^2\sin^2\rho}{a^2\sin^2\rho+b'^2\cos^2\rho}=\frac{a^2}{b'^2}, $$
therefore
$$ a^2b'^2\cos^2\rho+b'^4\sin^2\rho=a^4\sin^2\rho+a^2b'^2\cos^2\rho. $$
Simplifying, $a^4=b'^4$ and the result follows.
\end{proof}

\begin{lemma}\label{lem:final}
$\alpha=\pm i\alpha'$ and $\beta=\pm i \beta'$, with the same sign.
\end{lemma}
\begin{proof}
Recall \cref{eq:linear1},
\begin{equation*}
\big(\alpha \cos(\theta v_0) +\alpha' \sin(\theta v_0) \big)^{v_1}=\big(\beta\cos(\theta v_1)+\beta'\sin(\theta v_1)\big)^{v_0}.
\end{equation*}

If we consider the derivative on $\theta$ of \cref{eq:linear1} we obtain
\begin{equation}
\begin{aligned}
&v_0v_1\big(\alpha \cos(\theta v_0) +\alpha' \sin(\theta v_0) \big)^{v_1-1}\big(-\alpha \sin(\theta v_0) +\alpha' \cos(\theta v_0) \big)=\\
&v_0v_1\big(\beta\cos(\theta v_1)+\beta'\sin(\theta v_1)\big)^{v_0-1}\big(-\beta\sin(\theta v_1)+\beta'\cos(\theta v_1)\big).
\end{aligned}
\label{eq:linear2}
\end{equation}
Now we perform a trick in \cref{eq:linear2}: multiply by $\alpha \cos(\theta v_0) +\alpha' \sin(\theta v_0)$ and then simplify using \cref{eq:linear1}. This yields
\begin{equation}
\begin{aligned}
&\big(\beta\cos(\theta v_1)+\beta'\sin(\theta v_1)\big)\big(-\alpha \sin(\theta v_0) +\alpha' \cos(\theta v_0) \big)=\\
&\big(\alpha \cos(\theta v_0) +\alpha' \sin(\theta v_0) \big)\big(-\beta\sin(\theta v_1)+\beta'\cos(\theta v_1)\big).
\end{aligned}
\label{eq:lineargeneral}
\end{equation}
Evaluating \cref{eq:lineargeneral} at $\theta=0$ yields
\begin{equation}
\beta\alpha'=\alpha\beta'.
\label{eq:alphabeta2}
\end{equation}
Evaluating \cref{eq:lineargeneral} at $\theta=\frac{\pi}{2v_0}$ and using \cref{eq:alphabeta2}, after a small manipulation, it is easy to obtain
\begin{equation}
\frac{\beta\cos\left(\frac{\pi v_1}{2v_0}\right)+\beta'\sin\left(\frac{\pi v_1}{2v_0}\right)}{-\beta\sin\left(\frac{\pi v_1}{2v_0}\right)+\beta'\cos\left(\frac{\pi v_1}{2v_0}\right)}=\frac{\alpha}{\alpha'}=\frac{\beta}{\beta'}.
\label{eq:alphabetafinal}
\end{equation}
The result follows from \cref{eq:alphabetafinal,lem:orthogonal,lem:orthogonalcomplex}.
\end{proof}
   
   \subsection{Proof}
   We know from \cref{subsec:linear} that   
   $$ T_\infty\psi (x_R,x_I,y_R, y_I)=\left(\alpha x_R+\alpha'x_I,\beta y_R+\beta'y_I\right). $$
  This mapping is holomorphic if, by Cauchy--Riemann equations,  $\alpha=-i\alpha'$ and $\beta=-i\beta'$ and it is antiholomorphic if $\alpha=i\alpha'$ and $\beta=i\beta'$. By \cref{lem:final}, $T_\infty\psi$ is a holomorphic or antiholomorphic mapping. This is precisely what we mean by \textsl{rigidity}: the normal forms force the Taylor series of $\psi$ to be a holomorphic or antiholomorphic linear map.

\begin{theorem}
\label{thm:main}
Let $(C,0),(C',0)\subset(\CC^2,0)$ be two plane curve singularities, not necessarily irreducible, with equations $g(x,y)$ and $g'(x,y)$. Assume also that they contain a singular branch. Then, they are $\Ccal^\infty$-isomorphic if, and only if, they are $\CC$-isomorphic or $C$ is $\CC$-isomorphic to $\overline{C'}$, which is also a complex curve with equation $\overline{g'(\overline{x},\overline{y})}=0$.

 Furthermore, any diffeomorphism $\psi$ taking $C$ to $C'$, has a holomorphic or antiholomorphic Taylor series $T_\infty\psi$.
\end{theorem}
\begin{proof}
The same argument we have shown above works, and is far simpler, for monomial parametrizations $\gamma(t)=(t^{v_0},t^{v_1})$.

Assume that $\psi$ is the diffeomorphism that takes $C$ to $C'$, so
$$ g'\circ\psi\circ\gamma_i=0$$
for a parametrization $\gamma_i$ of any branch of $C$. In particular, $\psi$ takes a singular branch of $C$ to a singular branch of $C'$, so $T_\infty\psi$ is a holomorphic or antiholomorphic mapping by \cref{subsec:linear,lem:final}. It still holds that
$$ g'\circ T_\infty\psi\circ\gamma_i=0 $$
as formal power series, for any other branch of $C$. Hence, $T_\infty\psi$ also takes $C$ to $C'$. This proves the result.
\end{proof}

   With the Normal Forms Theorem (\cref{thm: normalform}) and \cref{thm:main} we can easily find pairs of curves that are $\Ccal^\infty$-isomorphic but not $\CC$-isomorphic or curves such that the $\Ccal^\infty$-class coincides with its $\CC$-class.
   
   \begin{example}\label{ex:moduli}
 The family of curves
\begin{equation}
\gamma_\varepsilon(t)=\bigg(t^4,t^9+t^{10}+\frac{19}{18} t^{11}+\varepsilon t^{15}\bigg)
\end{equation}
is, for any $\varepsilon\in \mathbb{C}$, a branch with semigroup $\Gamma = \langle 4,9\rangle$ with the property that the branches $\gamma_\varepsilon$ and $\gamma_{\varepsilon'}$ are analytically equivalent if and only if $\varepsilon=\varepsilon'$. This can be deduced from \cref{thm: normalform} (see also the 10-th row of \cite[Table 1]{Hefez2009} and \cite{Bayer2001}). This is therefore a parametric family of branches with different analytic classes for each value of $\varepsilon$. In particular, $\gamma_\varepsilon$ and $\gamma_{\overline{\varepsilon}}$ are always $\Ccal^\infty$-isomorphic but they are $\CC$-isomorphic if, and only if, $\varepsilon\in\RR$, in which case they parametrize the same set.
   \end{example}
   
 \section{Curves with only smooth branches}\label{s:smoothbranches}

\cref{thm:main} omits the case of several smooth branches. Let us consider a curve singularities $(C,0),(C',0)\subset(\CC^2,0)$ given by several lines. If $C$ and $C'$ have three lines, they are $\CC$-isomorphic, so we can assume that they have at least four lines. After long but straightforward computations of linear algebra, one can show that a $\Ccal^\infty$-isomorphism $\psi$ between $C$ and $C'$ with at least four lines induces a $\CC$-isomorphism between $C$ and $C'$ or $\overline{C'}$. Indeed, the differential $d\psi$ preserves three of the complex lines and it induces the $\CC$-isomorphism. We leave the details to the reader. This shows that the $\Ccal^\infty$-classification of \textit{tangent spaces} coincides with the $\CC$-classification modulo conjugation. 

However, in the case of only smooth branches, the Taylor series of a $\Ccal^\infty$-isomorphism could be far from being analytic. For example, one can take $\id\times f(y^R,y^I)$ as $\Ccal^\infty$-isomorphism taking the line $L=\left\{y=0\right\}$ to itself, where the only condition on $f$ is to be a diffeomorphism preserving the origin. One can construct a similar example for an arbitrary number of transverse lines and a radial diffeomorphism. Moreover, it is not even true that the $\Ccal^\infty$-classification and the $\CC$-classification coincide modulo conjugation when one considers tangent smooth branches, as the following example shows.

\begin{figure}[htbp]
      \centering
         \includegraphics[width=0.85\textwidth]{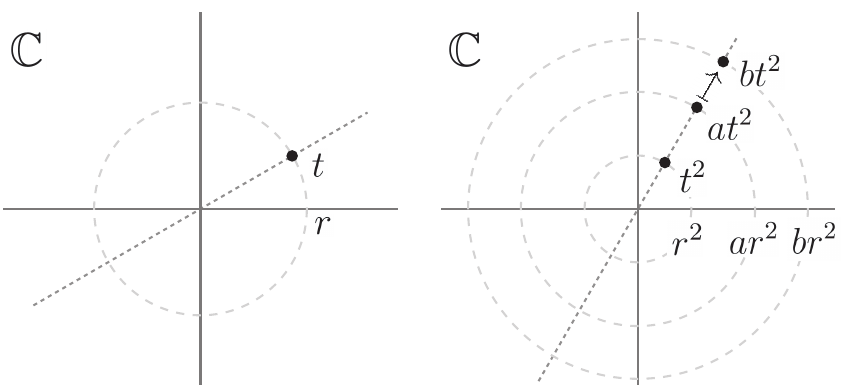}
      \caption{The diffeomorphism $\psi$ of \cref{eq:difeoexample}.}
      \label{fig:difeonew}
   \end{figure}

\begin{example}\label{ex:difeonew}
Consider the pair of curves with two branches (with $a,b\neq1$)
\begin{equation*}\begin{aligned}
(y-x^2)(y-ax^2)&=0\text{, and}\\
(y-x^2)(y-bx^2)&=0.
\end{aligned}\end{equation*}
It is shown in \cite[Section 5]{Hefez2015} (cf. \cite[Example 3]{Campillo1998}) that they are $\CC$-isomorphic if, and only if, $a=b$ or $a=\frac{1}{b}$. Assume they are not $\CC$-isomorphic and, for simplification, assume $a,b\in\RR$. It is clear that the branches of the curves are parametrized, respectively, by 
$$\begin{aligned}
t \mapsto(t,t^2),\ & t \mapsto(t,at^2) \text{ and }\\
t \mapsto(t,t^2),\ & t \mapsto(t,bt^2).
\end{aligned}$$

Consider the $\Ccal^\infty$-isomorphism $\psi:(\CC^2,0)\to(\CC^2,0)$ so that 
\begin{equation}\label{eq:difeoexample}
\psi(r_1 e^{i\theta_1},r_2e^{i\theta_2})=\big(r_1e^{i\theta_1}, \big(r_2 + \frac{b-1}{a-1}(r_2-r_1^2)\big)  e^{i\theta_2}\big),
\end{equation}
in polar coordinates (see \cref{fig:difeonew}). Writing $t=re^\theta$, $\psi$ fixes the branch $y-x^2=0$, since 
$$\psi(re^\theta,r^2 e^{i2\theta})=\big(re^\theta,\big(r^2+\frac{b-1}{a-1}(r^2-r^2) \big)e^{i2\theta}\big),$$
and sends $y-ax^2=0$ to $y-bx^2=0$, since
$$\psi(re^\theta,ar^2 e^{i2\theta})=\big(re^\theta,\big(r^2+\frac{b-1}{a-1}(ar^2-r^2) \big)e^{i2\theta}\big).$$
Hence, the curves are $\Ccal^\infty$-isomorphic, but they are not $\CC$-isomorphic. Observe that, when $a,b\in\RR$, they are also invariant by conjugation.
\end{example}

It is clear that more sophisticated examples with $\Ccal^\infty$-isomorphisms similar to \cref{eq:difeoexample} can be build from the ideas of \cref{ex:difeonew,fig:difeonew}. Hence, the $\Ccal^\infty$ classification of curves with only smooth branches is very coarse.

 \section{Hypersurfaces case}\label{sec:ephraim}
   
   We are now going to review Ephraim's work on \cite{Ephraim1973} and compare it with our own result. Recall that the main result of that paper is the following:
   
   \begin{theorem}[see {\cite[Theorem 4.2]{Ephraim1973}}]\label{thm:ephraim}
   Let $V$ be an irreducible germ of hypersurface at a point. Then, if $V$ is $\Ccal^\infty$-isomorphic to a germ of a complex set $W$, $V$ is $\CC$-isomorphic to $W$ or $\overline{W}$.
   \end{theorem}
   
   There are two reasons to review his work. The first one is because it is an interesting paper with some smart ideas that is not well-known enough in the singularity community. The second reason is that we can squeeze its proof to solve a big portion of the non-irreducible cases with similar ideas.
   \newline
   
   The first result he proves simplifies the problem of $\Ccal^\infty$-isomorphism to real analytic isomorphism.
   
   \begin{proposition}[see {\cite[Proposition 1.1]{Ephraim1973}}]\label{prop:ephranal}
   Let $V$ and $W$ be real analytic sets. If they are $\Ccal^\infty$-isomorphic then they are real analytic isomorphic.
   \end{proposition}

Ephraim's proof uses the fact that the rings of real analytic mappings of $V$ and $W$, $\Acal_V$ and $\Acal_W$, are isomorphic (which is equivalent to the statement) if their \textsl{completions}, $\widehat{\Acal_V}$ and $\widehat{\Acal_W}$, are isomorphic. This is a result of Artin (see \cite[Corollary 1.6]{Artin1968} and also the previous work of Hironaka and Rossi in \cite{Hironaka1964}). In general, for a local ring $R$, $\widehat{R}$ denotes its \textit{completion}, i.e., the (\textsl{projective}) limit of the (\textsl{projective}) system
      $$ \nicefrac{R}{\mfrak}\leftarrow\nicefrac{R}{\mfrak^2}\leftarrow\nicefrac{R}{\mfrak^3}\leftarrow\cdots $$
      
Despite not being explicit, it is not difficult to see along Ephraim's proof that, if $\psi$ is the diffeomorphism taking $V$ to $W$, then its Taylor series $T_\infty\psi$ induces the isomorphism between $\widehat{\Acal_V}$ and $\widehat{\Acal_W}$. Then, Artin's work implies that there are isomorphisms between $\Acal_V$ and $\Acal_W$ that coincide with $T_\infty\psi$ up to any arbitrary order. It is, however, unclear from this approach whether $T_\infty\psi$ is analytic itself. Therefore, since Ephraim's proof of \cref{thm:ephraim} starts with a real analytic isomorphism, his work is not constructive. This is the main difference with our proof of \cref{thm:ephraim} given for the plane curve case in \cref{thm:main}: we show the stronger statement that $T_\infty\psi$ itself is either holomorphic or antiholomorphic.
\newline

The next result provides an algebraic relation between $V$ and $W$ to prove \cref{thm:ephraim}.

\begin{theorem}[see {\cite[Theorem 2.1]{Ephraim1973}}]\label{thm:relations}
Let $(f_1,\dots,f_s)$ be the radical ideal of the germ of complex variety $(V,0)\subset(\CC^n,0)$ and $f_i^R,f_i^I\in\Acal_n$ the real and imaginary parts of $f_i$, respectively. The ideal of real analytic germs vanishing on $V$ is $(f_1^R,\dots,f_s^R,f_1^I,\dots,f_s^I)$ if, and only if, $V$ is irreducible. Equivalently, if $\overline{\bullet}$ denotes conjugation, the ideal of real analytic germs vanishing on $V$ is $(f_1 ,\dots,f_s ,\overline{f_1},\dots,\overline{f_s})$ if, and only if, $V$ is irreducible.
\end{theorem}

With this result, if $g_V,g_W$ are the reduced equations of $(V,0),(W,0)\subset(\CC^n,0)$ and $\psi$ is a real analytic isomorphism taking $V$ to $W$, then
\begin{equation}
\begin{aligned}
g_V\big(\psi^{-1}(\mathbf{z})\big)=&\eta_1(\mathbf{z})g_W(\mathbf{z})+\eta_2(\mathbf{z})\overline{g_W(z)}, and\\
g_W\big(\psi(\mathbf{z})\big)=&\xi_1(\mathbf{z})g_V(\mathbf{z})+\xi_2(\mathbf{z})\overline{g_V(z)}.
\end{aligned}
\label{eq:relations}
\end{equation}
At this point Ephraim uses an interesting technique to go from the real analytic category to the holomorphic category.
\begin{definition}[cf. {\cite[p. 18]{Ephraim1973}}]
Let $f\in\CC\left[\left[x_1,y_1,\dots,x_n,y_n\right]\right]$. Changing variables to $z_\bullet=x_\bullet+iy_\bullet,\overline{z_\bullet}=x_\bullet-iy_\bullet$ one can write $f=\sum_{I,J}a_{IJ}\mathbf{z}^I\overline{\mathbf{z}}^J$. The \textit{holomorphic Taylor series} of $f$ is 
$$ T_\mathbf{z}f=\sum_Ia_{I0}\mathbf{z}^I.$$
The definition for the \textit{antiholomorphic Taylor series} is analogous.
\end{definition}
It is not hard to see that:
\begin{itemize}
	\item $T_\mathbf{z}$ defines a ring homomorphism $T_\mathbf{z}:\CC\left[\left[x_1,y_1,\dots,x_n,y_n\right]\right]\to\CC\left[\left[z_1,\dots,z_n\right]\right]$,
   \item $T_\mathbf{z}$ restricts to $T_\mathbf{z}|:\Acal_{2n}\to\Ocal_n$, and
   \item $T_\mathbf{z}(f\circ g)=T_\mathbf{z}f\circ T_\mathbf{z}g$ (this is not mentioned in \cite{Ephraim1973}).
\end{itemize}
\begin{remark}
This is equivalent to form a Taylor-like series using the \textit{Wirtinger derivatives}:
$$ \frac{\partial}{\partial z}\coloneqq\frac{1}{2}\left(\frac{\partial}{\partial x}-i\frac{\partial}{\partial y}\right) $$
for the holomorphic Taylor series and
$$ \frac{\partial}{\partial \overline{z}}\coloneqq\frac{1}{2}\left(\frac{\partial}{\partial x}+i\frac{\partial}{\partial y}\right) $$
for the antiholomorphic Taylor series.
\end{remark}

Now, we can take the holomorphic Taylor series in \cref{eq:relations} and, by Cauchy-Riemann equations, we have that
\begin{equation}
\begin{aligned}
g_V\big(T_\mathbf{z}\psi^{-1}(\mathbf{z})\big)=&T_\mathbf{z}\eta_1(\mathbf{z})g_W(\mathbf{z}),\text{ and}\\
g_W\big(T_\mathbf{z}\psi(\mathbf{z})\big)=&T_\mathbf{z}\xi_1(\mathbf{z})g_V(\mathbf{z}).
\end{aligned}
\label{eq:relationsholo}
\end{equation}
At this point, if one shows that $T_\mathbf{z}\psi$ or $T_\mathbf{z}\psi^{-1}$ are isomorphisms, then $V$ and $W$ are $\CC$-isomorphic. In order to follow this direction one needs two extra conditions that are, as we will see later in \cref{ex:notholds,ex:notholds2}, necessary. 

The first condition is technical:
\begin{enumerate}
	\item[\textit{i)}] $\eta_1(0)\xi_1(0)\neq0$.
\end{enumerate}
 In this case $V$ and $W$ will be $\CC$-isomorphic. In a similar way, but arguing with the antiholomorphic Taylor series and conjugating \cref{eq:relations}, one can show that if
\begin{enumerate}
	\item[\textit{i')}] $\eta_2(0)\xi_2(0)\neq0$
\end{enumerate}  
then $V$ and $\overline{W}$ are $\CC$-isomorphic. This is proven in \cite[Lemma 4.1]{Ephraim1973}. Ephraim then shows in a simple argument that \textit{i)} or \textit{i')} hold.

The second condition is that $V$ or $W$ have zero-dimensional isosingular locus, the origin. Otherwise, $V\cong V'\times\CC^s, W\cong W'\times\CC^{s'}$ and one reasons with $V'$ and $W'\times\CC^{s-s'}$ (assume $s>s'$).

Under these two conditions, either the holomorphic or the antiholomorphic Taylor series of the $\Ccal^\omega$-isomorphism have linear part of maximum rank and they provide the isomorphism between $V$ and $W$ or $\overline{W}$ that concludes the proof (see \cite[Theorem 3.1]{Ephraim1973}). A large part of the article is dedicated to prove this technical fact. This finishes the proof of \cref{thm:ephraim}.

Despite the fact that Ephraim focuses his attention solely on the irreducible case, the ideas that are implemented to prove his main result, \cref{thm:ephraim}, can be adapted to extend the scope to some non-irreducible hypersurfaces.

\begin{definition}[see {\cite[Definition 3.2]{Ephraim1976}}]
An irreducible variety is \textit{non-factorable} or \textit{indecomposable} if it cannot be written as a product of varieties of lower dimension.
\end{definition}

\begin{proposition}
 Let $V$ be a (not necessarily irreducible) germ of hypersurface at a point. Assume that an irreducible component of $V$, say $V_0$, is non-factorable.
 Then, if $V$ is $\Ccal^\infty$-isomorphic to a germ of a complex set $W$, $V$ is $\CC$-isomorphic to $W$ or $\overline{W}$.
\end{proposition}
\begin{proof}
$V$ and $W$ are still $\Ccal^\omega$-isomorphic by some $\psi$, Artin's result \cite[Corollary 1.6]{Artin1968} holds in the non-irreducible case and so does \cref{prop:ephranal}. If one applies the argument shown above to prove \cref{thm:ephraim} to the component $V_0$, it follows that either the holomorphic or the antiholomorphic Taylor series of $\psi$ has linear part with maximum rank, since either condition \textit{i)} or \textit{i')} hold. Thus, \cref{eq:relationsholo}, or its antiholomorphic version, applied to each irreducible component of $V$ also shows that the same Taylor series provides the isomorphism taking all the components of $V$ to the components of $W$,or $\overline{W}$. Hence, $V$ is either $\CC$-isomorphic to $W$ or $\overline{W}$.
\end{proof}

The next example shows that the proof assumes necessary conditions, otherwise the statement is false. See also \cref{s:smoothbranches}.

\begin{example}\label{ex:notholds}
Ephraim in \cite[Example 2]{Ephraim1973} provides a family of normal surface singularities in $(\CC^3,0)$ that are not $\CC$-isomorphic to their conjugates. His family is a bit technical, but one can take the family given by four lines with different cross-ratios and exclude conjugates of this family (see also the family given in \cref{ex:moduli}). If one takes $V,V'$ two non-isomorphic elements of this family and considers 
$$\begin{aligned}
\mathbf{V}_1=&V\times\CC^2 \cup\CC^2\times V',\\ 
\mathbf{V}_2=&V\times\CC^2 \cup\CC^2\times\overline{V'},\\
\mathbf{V}_3=&\overline{V}\times\CC^2 \cup\CC^2\times V',\\
\mathbf{V}_4=&\overline{V}\times\CC^2 \cup\CC^2\times\overline{V'},
\end{aligned}$$ 
they are $\Ccal^\omega$-isomorphic but not $\CC$-isomorphic (this idea is also showed in \cite[Example 2]{Ephraim1973}).
For example, taking $\mathbf{V}_1$ and $\mathbf{V}_2$, the $\Ccal^\omega$-isomorphism is $\id_{\CC^2}\times\overline{\id_{\CC^2}}$. If one takes the (anti)holomorphic Taylor series of this map, it is $\id_{\CC^2}\times 0$ (or $0\times\id_{\CC^2}$). This shows that the argument, and the statement, are false for the non-irreducible case.
\end{example}

\section{Concluding remarks}\label{sec:conclude}

Let us consider the following example.

\begin{example}\label{ex:notholds2}
Consider any germ of irreducible hypersurface with zero-dimensional isosingular locus $(V,0)\subset(\CC^n,0)$. Now take $V\times\CC$. It is $\Ccal^\infty$-isomorphic to itself if one considers the map $\psi\coloneqq\id_{\CC^n}\times f$, where $f:(\CC,0)\to(\CC,0)$ is a (real) diffeomorphism with a Taylor series that does not converge at any point. These functions certainly exist, one can take examples of such functions from the real setting (see \cite{Kim2000}). Then $T_\infty\psi$ is obviously not analytic. However, $T_K\psi$ for any $K\geq1$ is analytic and a $\Ccal^\omega$-isomorphism of $V\times\CC$ to itself (not necessarily a $\CC$-isomorphism either).
\end{example}

In view of \cref{ex:notholds2} and our main result \cref{thm:main}, it is natural to ask the following question.

\begin{question}
Let $V$ be an irreducible germ of hypersurface at a point with zero-dimensional isosingular locus. Then, if $\psi$ is a diffeomorphism taking $V$ to a germ of a complex set $W$, the Taylor series of $\psi$ is either holomorphic or antiholomorphic.
\end{question}

This holds for plane curve singularities by \cref{thm:main} and it is not true for positive-dimensional isosingular locus by \cref{ex:notholds2}.

\begin{example}[see {\cite[Example 1]{Ephraim1973}}]\label{ex:c1notan}
Ephraim showed that the curves given by the parametrizations
$$
\begin{aligned}
\gamma_1(t)&=(t^3,t^7)\text{ and}\\
\gamma_2(t)&=(t^3,t^7+t^8)
\end{aligned} $$
are $\Ccal^1$-isomorphic, but not $\CC$-isomorphic. It is easy to see that they are topologically equivalent. In contrast, two sets of four lines with different cross-ratios are an example of plane curves that are topologically equivalent but not $\CC$, $\Ccal^1$, $\Ccal^\omega$ nor $\Ccal^\infty$-isomorphic.
\end{example}

As far as we know, $\Ccal^1$-equivalence is not closed, and it is the last natural classification one can give that remains open. Examples such as the one given above motivate this question, which could seem counterintuitive at first.

\begin{question}\label{qu:c1}
The $\Ccal^1$-equivalence coincides with the topological equivalence for irreducible plane curve singularities.
\end{question}

   \section*{Acknowledgements}

We are specially thankful to Marco Marengon, who motivated us to answer the question about smooth classification of curves and suffered many of the partial technical results in the process of making this work. We thank Marcelo Escudeiro Hernandes for his useful comments and encouragement. We also want to thank Edson Sampaio for pointing out to us Ephraim's work, unfortunately few days after finishing our proof.
   
   \section*{Data availability}
   This work does not use any set of data.
   \section*{Conflict of interest statement}
   The authors state that there is no conflict of interest in the making and publication of this work.
\bibliographystyle{C:/Users/rgc19/AppData/Local/Programs/MiKTeX/bibtex/bst/bibtex/myalpha.bst}
\bibliography{C:/Users/rgc19/Documents/Textos-Matematicas/Bibliografias/bibtex/bib/mybibs/FullBib.bib} 
\end{document}